\documentclass{amsart}

\usepackage[latin1]{inputenc}
\usepackage{amsmath,amssymb,graphicx,setspace,verbatim}
\usepackage{color}
\usepackage{comment}
\usepackage{appendix}
\theoremstyle{comment}
\usepackage[color,matrix,arrow]{xy}

\newtheorem*{mcomment}{\color{cyan}{Comment}}

\newcommand{\Sym}{\mathrm{S}} 

\DeclareMathOperator{\Sp}{Sp}
\DeclareMathOperator{\SL}{SL}

\DeclareMathOperator{\Orth}{O}
\DeclareMathOperator{\Alt}{Alt}

\newcommand{\la}{\langle}
\newcommand{\ra}{\rangle}
\newcommand{\Fk}{\Bbb{F}_{2^k}}

\newcommand{\Fp}{\Bbb{F}_p}
\renewcommand{\Sym}{{\rm Sym}}
\renewcommand{\leq}{\leqslant}
\renewcommand{\geq}{\geqslant}
\renewcommand{\bar}{\overline}

\input xy
\xyoption{all}

\newtheorem{theorem}{Theorem}[section]

\newtheorem{lemma}[theorem]{Lemma}
\newtheorem{corollary}[theorem]{Corollary}
\newtheorem{proposition}[theorem]{Proposition}
\theoremstyle{definition}

\newtheorem{conjecture}[theorem]{Conjecture}

\begin{document}

\title{Rank reduction of string C-group representations}

\author{Peter A. Brooksbank}
\address{Peter A. Brooksbank, 
	Department of Mathematics,
	Bucknell University,
	Lewisburg, PA 17837,
	USA
}
\email{pbrooksb@bucknell.edu}
\author{Dimitri Leemans}
\address{Dimitri Leemans, Universit\'{e} Libre de Bruxelles,
D\'{e}partement de Math\'{e}matique,
C.P.216 Alg\`ebre et Combinatoire,
Bld du Triomphe, 1050 Bruxelles,
Belgium}
\email{dleemans@ulb.ac.be}

\thanks{This work was partially supported by a grant from the 
Simons Foundation (\#281435 to Peter Brooksbank), and by
the Hausdorff Research Institute for Mathematics.}

\maketitle
\begin{abstract}
We show that a rank reduction technique for string C-group representations first 
used in~\cite{fl} for the symmetric groups generalizes to arbitrary settings. The 
technique permits us, among other things, to prove that orthogonal groups defined 
on $d$-dimensional modules over fields of even order greater than 2 possess string 
C-group representations of all ranks $3\leq n\leq d$. The broad applicability of the 
rank reduction technique provides fresh impetus to construct, for suitable families of 
groups, string C-groups of highest possible rank. It also suggests that the alternating 
group $\Alt(11)$---the only known group having `rank gaps'---is perhaps
 more unusual than previously thought.
\end{abstract}

\noindent \textbf{Keywords:} abstract regular polytope, string C-group, Coxeter group.

\noindent \textbf{2000 Math Subj. Class:} 52B11, 20D06.

\section{Introduction}
In a recent joint work of the second author and Fernandes~\cite{fl}, a certain ``rank reduction" 
technique was used to show that the symmetric group $\Sym(m)$ acts as the 
group of automorphisms of an abstract regular polytope of rank $n$ for every 
$3\leq n\leq m$. The technique was used again by those same to authors in~\cite{fl2} 
to prove a similar result for the alternating group $\Alt(m)$. Both applications of rank 
reduction seemed at the time to depend crucially on their particular setting, namely 
$\Sym(m)$ or $\Alt(m)$ acting on the natural permutation domain $\{1,.\ldots,m\}$. In this 
paper we show, to the contrary, that it is a substantially more general technique.

Abstract regular polytopes have an equivalent formulation in terms of quotients
of Coxeter groups with string diagrams, and it is helpful to frame our discussion in these terms.
We say that  $(G;\{\rho_0,\ldots, \rho_{n-1}\})$ is a 
{\em string group generated by involutions}---or {\em sggi} for short---if
 $G=\la \rho_0,\ldots,\rho_{n-1}\ra$, each $\rho_i$ is an involution, and the 
 sequence $\rho_0,\ldots,\rho_{n-1}$ satisfies the {\em commuting property}
 \begin{align}
 \label{eq:comm-prop}
 \forall i,j\in\{0,\ldots, n-1\} & &  \;|i-j|>1\Rightarrow (\rho_i\rho_j)^2=1.
 \end{align}
 If an sggi $(G;\{\rho_0,\ldots, \rho_{n-1}\})$ additionally satisfies the {\em intersection property}
 \begin{align}
 \label{eq:int-prop}
 \forall I,J \subseteq \{0,\ldots,n-1\} & &  \langle \rho_i \mid i \in I\rangle \cap \langle \rho_j \mid j 
 \in J\rangle = \langle \rho_k \mid k \in I\cap J\rangle
 \end{align}
 then it is a {\em string C-group}, 
 and $n$ is its {\em rank}.
 If $|\rho_i\rho_{i+1}|>2$ for all $0\leq i\leq n-2$ then 
 $(G;\{\rho_0,\ldots, \rho_{n-1}\})$ 
 is {\em irreducible};
 unless $G$ is directly decomposable, the string C-group must be irreducible.
 Our main result is the following.

%%%%
\begin{theorem}[Rank Reduction]
\label{thm:main}
Let $(G;\{\rho_0,\ldots, \rho_{n-1}\})$ be an irreducible string C-group of rank $n\geq 4$.
If $\rho_0\in \langle \rho_0\rho_2,\rho_3\rangle$, then 
$(G;\{\rho_1, \rho_0\rho_2, \rho_3, \ldots, \rho_{n-1}\})$ is a string C-group of rank $n-1$.
\end{theorem}

The condition $\rho_0\in \langle \rho_0\rho_2,\rho_3\rangle$ is an easy one to verify, making the
Rank Reduction Theorem a powerful tool in the search for new polytopes.
For example, suppose that $\rho_2\rho_3$ has odd order $2k+1$. Then 
\[
((\rho_0\rho_2)\rho_3)^{2k+1}=(\rho_0(\rho_2\rho_3))^{2k+1}=\rho_0\in \la \rho_0\rho_2,\rho_3\ra,
\]
so we obtain the following immediate and useful consequence of Theorem~\ref{thm:main}.

%%%%
\begin{corollary}
\label{coro:odd-order-single}
Let $(G;\{\rho_0,\ldots, \rho_{n-1}\})$ be an irreducible string C-group of rank $n\geq 4$.
If $\rho_2\rho_3$ has odd order, then 
$(G;\{\rho_1, \rho_0\rho_2, \rho_3, \ldots, \rho_{n-1}\})$ is a string C-group of rank $n-1$.
\end{corollary}
The integer sequence $[p_1,\ldots,p_{n-1}]$, where $p_i$ is the order of $\rho_{i-1}\rho_i$ 
is called the {\em Schl\"{a}fli type} of the string C-group $(G;\{\rho_0,\ldots, \rho_{n-1}\})$. So,
Corollary~\ref{coro:odd-order-single} tells us to look at the third integer $p_3$ 
in the Schl\"{a}fli type to see if we can apply the rank reduction mechanism once. 
Suppose we can, and suppose that $n\geq4$. Then
we obtain a new string C-group $(G;\{\rho_1, \rho_0\rho_2, \rho_3, \ldots, \rho_{n-1}\})$
with Schl\"{a}fli type $[q_1,\ldots,q_{n-2}]$, where $q_1$ is the order of $\rho_1\rho_0\rho_2$,
$q_2$ is the order of $\rho_0\rho_2\rho_3$, and $q_i=p_{i+1}$ for $3\leq i\leq n-2$.
Thus, if $n>4$ and $p_4$ is also odd, we can repeat the rank reduction to obtain
a string C-group of rank $n-2$. Iterating, we obtain the following result.

%%%%
\begin{corollary}
\label{coro:odd-order-sequence}
Let $(G;\{\rho_0,\ldots, \rho_{n-1}\})$ be an irreducible string C-group of rank $n\geq 4$.
Let $[p_1,\ldots, p_{n-1}]$ be its Schl\"afli type, and put
\[
t=\max\{j\in\{0,\ldots,n-3\}\colon \forall i\in\{0,\ldots,j\},\; p_{2+i}\;\mbox{is odd}\}.
\]
Then $G$ is a string C-group of rank $n-i$ for each $i\in\{0,\ldots,t\}$.
\end{corollary}

Section~\ref{sec:proof} is devoted to the proof of the Rank Reduction Theorem. 
While the condition $\rho_0\in\la \rho_0\rho_2,\rho_3\ra$ is 
convenient to ensure the success of the process, it is not essential;
we illustrate this in Section~\ref{sec:proof}. In Section~\ref{sec:sym-alt}, 
we show that some of the striking results proved in~\cite{fl} and~\cite{fl2} for the groups 
$\Sym(m)$ and $\Alt(m)$ are, in fact, immediate consequences of 
Corollary~\ref{coro:odd-order-sequence}. Indeed, the potential to iterate rank reduction
on string C-groups via Corollary~\ref{coro:odd-order-sequence} impresses on us the importance
of obtaining {\em \underline{some}} general high rank construction for suitable families of groups. 
To emphasize that point, in Section~\ref{sec:sp-orth}
we examine recent constructions of high rank polytopes for symplectic and orthogonal groups through the
lens of rank reduction, proving in particular the following new result.

%%%%%
\begin{theorem}
\label{thm:sp}
Let $k\geq 2$ and $m\geq 2$ be integers.
\begin{enumerate}
\item[(a)]  The symplectic group $\Sp(2m,\Fk)$ is a string C-group of rank $n$ for each $3\leq n\leq 2m+1$.
 \item[(b)]  The orthogonal groups $\Orth^{+}(2m,\Fk)$ and $\Orth^-(2m,\Fk)$ are
 string C-groups of rank $n$ for each $3\leq n\leq 2m$.
 \end{enumerate}
\end{theorem}

%%%%%%%
%%%%%%%
\section{The Rank Reduction Theorem}
\label{sec:proof}
Let $(G;\,\{\rho_0,\ldots,\rho_{n-1}\})$ be an sggi of rank $n$. 
If $G$ is a string C-group,
and $I\subseteq \{0,\ldots,n-1\}$, it readily follows that 
$\la \rho_i\colon i\in I\ra$ is also a string C-group on
its defining generating sequence. The following 
result, proved in~\cite[Proposition 2E16]{arp}, facilitates an inductive 
approach to verifying string  C-group representations.

%%%
\begin{proposition}
\label{arp}
Let $(G;\,\{\rho_0,\ldots,\rho_{n-1}\})$ be an sggi.
If the subgroups $\la \rho_0,\ldots,\rho_{n-2}\ra$
and $\la \rho_1,\ldots,\rho_{n-1}\ra$ are both string C-groups, and
\[
\la \rho_0,\ldots,\rho_{n-2}\ra \cap \la \rho_1,\ldots,\rho_{n-1}\ra =
\la \rho_1,\ldots,\rho_{n-2}\ra,
\]
then $(G;\,\{\rho_0,\ldots,\rho_{n-1}\})$ is a string $C$-group.
\end{proposition}

The next elementary result follows from~\cite[Classification Theorem 1.2]{CL2015}.

%%%
\begin{lemma}
\label{lem:dihedral}
The dihedral group of order $2k\geq 6$ is an irreducible string C-group of rank $n$
if, and only if, $n=2$.
\end{lemma}

\noindent {\bf Proof of Theorem~\ref{thm:main}.}~
From the hypothesis $\rho_0\in\la \rho_0\rho_2,\rho_3\ra$ it is immediate that
rank reduction on the generators of $G$ does not yield a proper subgroup. 
It remains to show that the new generating sequence is again a string C-group
representation of $G$. For this we induct on the rank $n\geq 4$.

In the base case we have an irreducible string C-group $(G;\,\{\rho_0,\rho_1,\rho_2,\rho_3\})$ 
of rank 4, and must show that $(G;\,\{\rho_1,\rho_0\rho_2,\rho_3\})$ is a string C-group. Evidently,
both $\la \rho_1,\rho_0\rho_2 \ra$ and $\la \rho_0\rho_2,\rho_3\ra$ are string C-groups,
so it suffices to show that the intersection of these two dihedral groups is the cyclic group 
$\la \rho_0\rho_2\ra$. Note,
\[
\la \rho_1,\rho_0\rho_2 \ra \cap \la \rho_0\rho_2,\rho_3\ra
\leq
\la \rho_0,\rho_1,\rho_2 \ra \cap \la \rho_0,\rho_2,\rho_3 \ra
=
\la \rho_0,\rho_2\ra,
\]
since $G$ is a string C-group. Thus, if $\la \rho_1,\rho_0\rho_2 \ra \cap \la \rho_0\rho_2,\rho_3\ra$
properly contains $\la \rho_0\rho_2\ra$, it must contain $\rho_0$. But then
$\la \rho_0,\rho_1,\rho_2\ra=\la \rho_1,\rho_0\rho_2\ra$ is dihedral, contrary to Lemma~\ref{lem:dihedral}.

Next, suppose $(G;\,\{\rho_0,\ldots,\rho_{n-1}\})$ is an irreducible string C-group of rank
$n>4$, with $\rho_0\in\la \rho_0\rho_2,\rho_3\ra$.
Assume the result holds for string C-groups of smaller ranks.
In particular, since 
$H:=\la \rho_0,\ldots,\rho_{n-2}\ra$ 
is an irreducible string C-subgroup of $G$ of rank $n-1$,
by induction it follows that
$(H;\;\{ \rho_1,\rho_0\rho_2,\rho_3,\ldots,\rho_{n-2} \})$
is a string C-group (of rank $n-2$). Further, since $\rho_0\in\la \rho_0\rho_2,\rho_3\ra$, it follows that
\[
K:=\la \rho_0\rho_2,\rho_3,\ldots,\rho_{n-1} \ra=\la \rho_0,\rho_2,\ldots,\rho_{n-1}\ra.
\]
As $(K;\;\{ \rho_0,\rho_2,\rho_3,\ldots,\rho_{n-1} \})$ is a string C-group,
an easy induction shows that $(K;\;\{ \rho_0\rho_2,\rho_3,\ldots,\rho_{n-1} \})$
is also a string C-group. 
Finally,
\[
\begin{array}{rcl}
H\cap K & = & \la \rho_1,\rho_0\rho_2,\rho_3,\ldots,\rho_{n-2} \ra \cap \la \rho_0\rho_2,\rho_3,\ldots,\rho_{n-1} \ra \\
              & = &  \la \rho_0,\ldots,\rho_{n-2}\ra \cap \la \rho_0,\rho_2,\ldots,\rho_{n-1}\ra \\
              & = & \la \rho_0,\rho_2,\rho_3,\ldots,\rho_{n-2} \ra \\
              & = & \la \rho_0\rho_2,\rho_3,\ldots,\rho_{n-2}\ra,
\end{array}
\]
and it now follows from Proposition~\ref{arp} that $(G;\,\{\rho_1,\rho_0\rho_2,\rho_3,\ldots,\rho_{n-1}\})$ is a string C-group,
as required.
$\Box$
\bigskip

\noindent {\bf The \boldmath$\rho_0\in\la \rho_0\rho_2,\rho_3\ra$ criterion.}~
The condition $\rho_0\in\la \rho_0\rho_2,\rho_3\ra$, while easy both to state and verify,
is by no means essential for the rank reduction trick to work. We conclude this section with an example
that can readily be checked on a computer.
Consider the group $G$ preserving the symmetric form
\[
F=
\left[ \begin{array}{rrrr}
 1 & 1 & 0 & 0 \\
 1 & 2 & 1 & 0 \\
 0 & 1 & 1 & 2 \\
 0 & 0 & 2 & 1
\end{array}
\right]
\]
defined over $\Bbb{F}_3$. The Witt index of $F$ is 1, so $G=\Orth^-(4,\Bbb{F}_3)\cong\SL(2,\Bbb{F}_9)\colon 2$.
For $v\in V=\Bbb{F}_3^{\,4}$ nonsingular, let $\tau(v)$ denote the reflection in the 1-space 
$\la v\ra$ relative to the form $F$. Let $v_0,v_1,v_2,v_3$ denote the standard basis of $V$ 
relative to which $F$ is written and,
for $i,j\in\{0,1\}$, put $\rho_{2i+j}:=(-1)^j\tau(v_{2i+j})$. As matrices,
\[
\rho_0=\left[ \begin{smallmatrix} 2 & 0 & 0 & 0 \\ 1 & 1 & 0 & 0 \\ 0 & 0 & 1 & 0 \\ 0 & 0 & 0 & 1 \end{smallmatrix}  \right],~~~
\rho_1=\left[ \begin{smallmatrix} 2 & 1 & 0 & 0 \\ 0 & 1 & 0 & 0 \\ 0 & 1 & 2 & 0 \\ 0 & 0 & 0 & 2 \end{smallmatrix}  \right],~~~
\rho_2=\left[ \begin{smallmatrix} 1 & 0 & 0 & 0 \\ 0 & 1 & 1 & 0 \\ 0 & 0 & 2 & 0 \\ 0 & 0 & 2 & 1 \end{smallmatrix}  \right],~~~
\rho_3=\left[ \begin{smallmatrix} 2 & 0 & 0 & 0 \\ 0 & 2 & 0 & 0 \\ 0 & 0 & 2 & 1 \\ 0 & 0 & 0 & 1 \end{smallmatrix}  \right].
\]
Then $(G;\,\{\rho_0,\rho_1,\rho_2,\rho_3\})$ is a string C-group with Schl\"{a}fli type $[4,4,6]$. 
Applying rank reduction, $(G;\,\{\rho_1,\rho_0\rho_2,\rho_3\})$
is also a string C-group of Schl\"{a}fli type $[6,6]$. Finally,
$\rho_0$ is not in the dihedral group $\la\rho_0\rho_2,\rho_3\ra$.

%%%%%%%%%%%%%%%%%%%%%%
\section{Symmetric and alternating groups}
\label{sec:sym-alt}
The $(m-1)$-simplex of $\Sym(m)$ has Schl\"{a}fli type consisting entirely of 3's.
Hence, applying Corollary~\ref{coro:odd-order-sequence} gives a 
family of string C-group representations for $\Sym(m)$ of every rank from 
$m-1$ down to 3. This was already proved in~\cite[Theorem 3]{fl}.

Recently, Fernandes and Leemans~\cite{fl2} showed for $m\geq 12$ that $\Alt(m)$ has string 
C-group representations of every rank $\lfloor (m-1)/2\rfloor$ down to 3.
However, only in the case $m\equiv 1\;({\rm mod}\;4)$ 
can one apply Corollary~\ref{coro:odd-order-sequence} to
 the string C-group representation of (highest) rank 
$ \lfloor (m-1)/2\rfloor$ to produce representations of all ranks down to 3.
We remark, though, that for $n$ even~\cite[Theorem 4.1]{fl2} is a consequence of 
Corollary~\ref{coro:odd-order-sequence}.

These findings raise the question, for a family of groups, of whether one can find
a highest rank string C-group representation that, via rank reduction, generates 
all other permissible ranks of string C-group representations. We have an affirmative 
answer for the symmetric groups and for subfamilies of alternating groups.

At the time of writing, the only group that has gaps in its set of possible ranks of 
string C-group representations is the alternating group $\Alt(11)$. This group has 
representations of ranks 3 and 6 but none of ranks 4 or 5.
Figure~\ref{A11} gives the Schl\"{a}fli types and 
CPR graphs\footnote{A CPR graph is a special type of Schreier coset graph that encodes 
the defining involutions of a string C-subgroup of $\Sym(m)$ as a labelled graph
on the points of its domain $\{1,\ldots,m\}$.}, extracted from~\cite{flm}, of 
the three pairwise non-isomorphic rank 6 string C-group representations of $\Alt(11)$. 
The reader will immediately notice that third integer (from the left or right!) in the Schl\"{a}fli
type is 6, so one cannot hope to use Corollary~\ref{coro:odd-order-single}.
Further, one can readily see from the graphs in Figure~\ref{A11} that applying our rank reduction 
technique (again, from the left or from the right) will produce  
graphs that are not connected. This means that our rank reduction mechanism always
generates proper (intransitive) subgroups of $\Alt(11)$.

\begin{figure}
$$\begin{tabular}{| c | c |}
\hline 
Type & CPR graph  \\ \hline \hline
\{5,3,6,3,5\} &  \xymatrix@-1.1pc{*+[o][F]{} \ar@{=}[r]^2_0 & *+[o][F]{} \ar@{-}[r]^1 & *+[o][F]{} \ar@{-}[r]^0 & *+[o][F]{} \ar@{-}[r]^1 & *+[o][F]{} \ar@{-}[r]^2 & *+[o][F]{} \ar@{-}[r]^3 &*+[o][F]{} \ar@{-}[r]^{4} & *+[o][F]{} \ar@{-}[r]^{5} & *+[o][F]{} \ar@{-}[r]^{4} & *+[o][F]{} \ar@{=}[r]^{5}_3 & *+[o][F]{}  } \\ 
\hline
\{5,5,6,3,5\} &  \xymatrix@-1.1pc{*+[o][F]{} \ar@{-}[r]^0 & *+[o][F]{} \ar@{-}[r]^1 & *+[o][F]{} \ar@{=}[r]^0_2 & *+[o][F]{} \ar@{-}[r]^1 & *+[o][F]{} \ar@{-}[r]^2 & *+[o][F]{} \ar@{-}[r]^3 &*+[o][F]{} \ar@{-}[r]^{4} & *+[o][F]{} \ar@{-}[r]^{5} & *+[o][F]{} \ar@{-}[r]^{4} & *+[o][F]{} \ar@{=}[r]^{5}_3 & *+[o][F]{} } \\ 
\hline
\{5,5,6,5,5\} &  \xymatrix@-1.1pc{*+[o][F]{} \ar@{-}[r]^0 & *+[o][F]{} \ar@{-}[r]^1 & *+[o][F]{} \ar@{=}[r]^0_2 & *+[o][F]{} \ar@{-}[r]^1 & *+[o][F]{} \ar@{-}[r]^2 & *+[o][F]{} \ar@{-}[r]^3 &*+[o][F]{} \ar@{-}[r]^{4} & *+[o][F]{} \ar@{=}[r]^{5}_3 & *+[o][F]{} \ar@{-}[r]^{4} & *+[o][F]{} \ar@{-}[r]^{5} & *+[o][F]{}} \\ 
\hline
\end{tabular}$$

\caption{CPR graphs of rank six regular polytopes for $\Alt(11)$.}\label{A11}
\end{figure}

%%%%%%%
%%%%%%%
\section{Orthogonal and symplectic groups}
\label{sec:sp-orth}
The previous section provided a striking illustration of the power of rank reduction 
in situations where we can get our hands on at least one string C-group representation 
of high rank. Unfortunately, other than the symmetric and alternating groups, very few 
families of groups are known to have such representations. In this final section, however, 
we revisit two constructions for orthogonal and symplectic groups, applying our rank 
reduction technique to obtain new results.
\bigskip

%%%%%%%
\noindent {\bf Modular reduction of string Crystallographic groups.}
In a series of three papers, Monson and Schulte~\cite{MS1,MS2,MS3} conducted a comprehensive 
investigation of string C-groups that arise under modular reduction of crystallographic Coxeter groups.
Their setting is as follows. Let $\Gamma=\la \rho_0,\ldots,\rho_{n-1}\ra$ be an abstract Coxeter group
with string diagram and Schl\"{a}fli type $[p_1,\ldots,p_{n-1}]$. Suppose $\Gamma$ has a faithful
representation $\Gamma\to{\rm GL}(n,\Bbb{Z})$ as a group of reflections, and let $G=\la r_0,\ldots,r_{n-1}\ra$
be its image in ${\rm GL}(n,\Bbb{Z})$. For each odd prime $p$, one can consider the group 
 $\bar{G}=\la \bar{r}_0,\ldots,\bar{r}_{n-1}\ra\leq{\rm GL}(n,\Fp)$ obtained by reducing matrix entries modulo $p$.
 Evidently, $\bar{G}$ remains an sggi. Further, although $\bar{G}$ does not automatically inherit the intersection property 
 from its parent group $G$, Monson and Schulte show that in many cases it does.
 
 Of particular interest to us are the constructions of high rank. In~\cite[Section 3]{MS3} the authors study the family
 of so-called `3--infinity' groups, namely those having Schl\"{a}fli type $[3^k,\infty^l,3^m]$. Setting $m=0$, they show in
 particular that for any prime $p\geq 5$, and any positive integers $k,l$,
 the reduction $\bar{G}$ modulo $p$ of the string crystallographic Coxeter group
 $G$ of Schl\"{a}fli type $[3^k,\infty^l]$ is again a string C-group.
 Furthermore, $\bar{G}$ is a maximal subgroup of $\Orth(k+l+1,\Fp)$
 and has Schl\"{a}fli type $[3^k,p^l]$. The residue class of $p$ modulo 4 determines which 
 maximal subgroup we generate---the generating reflections all have the same spinor norm
 relative to the quadratic form preserved by $\bar{G}$, and that norm varies according to $p$.
 (Note, we ignore the reduction modulo 3 since that just gives a string C-group representation of 
 $\Sym(k+l+1)$ as permutation matrices.)
 
 A direct application of Corollary~\ref{coro:odd-order-sequence} now yields the following result:
 
 %%%%
 \begin{theorem}
 For each $d\geq 3$ and each odd prime $p\geq 5$, there is a subgroup $M$ of $\Orth(d,\Fp)$ of index $2$
 such that $M$ is a string C-group of rank $n$ for each $3\leq n\leq d$.
 \end{theorem}
 \medskip
 
 \noindent {\bf Proof of Theorem~\ref{thm:sp}.}
 More recently, the authors showed in joint work with Ferrara that, provided $k\geq 2$, the groups 
 $\Orth(2m,\Fk)$ and $\Sp(2m,\Fk)$  also have string C-group representations of high 
 rank~\cite[Corollary 1.5]{Brooksbank:2018}: rank $2m$ for orthogonal groups, and rank $2m+1$ for
 symplectic groups, since $\Sp(2m,\Fk)\cong\Orth(2m+1,\Fk)$. The approach---in some sense dual to that of 
 Monson and Schulte---was to construct a quadratic form $\varphi$ as a matrix relative to a basis
 of nonsingular vectors in such a way that the sequence of {\em symmetries} (generalizations of reflections) determined by
 these vectors gives a generating sequence for the desired orthogonal group as a string C-group. 
 (The example given at the end of Section~\ref{sec:proof} was based on this construction.)
 It was further shown,
 in even dimension, how to control the Witt index of the orthogonal space determined by $\varphi$. It is a natural
 consequence of the construction that the product of successive symmetries in the generating sequence is
 always an element of order $2^k+1$~\cite[Section 5]{Brooksbank:2018}. 
 In particular, the Schl\"{a}fli types of the string C-group representations arising from this
 construction are sequences of odd integers. Hence, Theorem~\ref{thm:sp} is an immediate consequence of
 Corollary~\ref{coro:odd-order-sequence} and~\cite[Corollary 1.5]{Brooksbank:2018}.

\section{Final thoughts}
Our rank reduction theorem shows the real interest in trying to find, for a given infinite family of groups, 
the string C-group representations of highest possible ranks. Indeed, given a `highest rank' string C-group 
representation for a group $G$, we can attempt to use the rank reduction technique to produce
new string C-group representations of lower ranks.

We have seen how this works almost effortlessly for the family of symmetric groups $\Sym(m)$,
orthogonal groups $\Orth^{\pm}(2m,\Fk)$, and symplectic groups $\Sp(2m,\Fk)$.
We have also seen how, with substantially more effort, it can be used to fill in the `rank gaps' in the 
alternating groups $\Alt(m)$.
The single exception in this regard is the group $\Alt(11)$, which was already identified 
in~\cite{fl2} as being special. In view of how broadly successful the rank reduction technique 
appears to be (albeit from our somewhat limited experience using it) the group $\Alt(11)$ strikes
us now as an anomaly, and prompts us to conclude with the following conjecture. 

\begin{conjecture}
The group $\Alt(11)$ is the only finite simple group whose set of ranks of string C-group 
representations is not an interval in the set of integers.
\end{conjecture}

\section{Acknowledgements}
We thank an anonymous referee for helpful comments and suggestions.

%-----------------------------------------------------------------------------------------------------------------------
%-----------------------------------------------------------------------------------------------------------------------
\bibliographystyle{amsplain}

\end{document}